\documentclass[12pt]{article}
\usepackage{amstex}
\usepackage{epsf}
\usepackage[all]{xy}
\usepackage[english,french]{babel}

\input xy
\xyoption{all}

\newtheorem{theo}{Theorem}[section]
\newtheorem{prop}[theo]{Proposition}
\newtheorem{cor}[theo]{Corollary}
\newtheorem{dfn}[theo]{Definition}

\newtheorem{rmk}[theo]{Remark}
\newtheorem{lemma}[theo]{Lemma}

\newenvironment{proof}{{\flushleft \bf Proof: \,} }

\newenvironment{ex}{{\flushleft \bf Example: \, }}

\def\ot{\otimes}
\def\gl{\Gamma_{\! \! \Lambda}}
\def\gln{\Gamma_{\! \! \Lambda _n}}

\def\lan{\Lambda _n}
\def\m{\frak{m}}

\def\uq{$u_q^+(\frak{sl_2})$}

\advance\oddsidemargin by- 2.1cm
\advance\topmargin by- 1.2cm 
\numberwithin{equation}{section}

\title{\textsc{Cyclic homology of the Taft algebras and of their Auslander algebras}}

\author{Rachel Taillefer\thanks{Laboratoire G.T.A., D\'epartement de Math\'ematiques CC 51, Universit\'e Montpellier II, 34095 Montpellier Cedex 5. 
email: taillefr\at math.univ-montp2.fr}}

\date{}

\begin{document}

\maketitle


\selectlanguage{english}
\begin{abstract} In this paper, we compute  the cyclic homology of the Taft algebras and of their Auslander algebras. Given a Hopf algebra $\Lambda,$ the Grothendieck groups of projective  $\Lambda -$modules and of all   $\Lambda -$modules are endowed with a  ring structure, which in the case of the Taft algebras is commutative (\cite{C2}, \cite{G}). We also describe the first Chern character for these algebras.
\end{abstract}


\section{Introduction}

The object of this paper is to compute the cyclic homology and the Chern characters of the Taft algebras $\lan$ and of their Auslander algebras $\gln,$ in order to study a possible influence of  the Hopf algebra structure of $\lan$ on them.

Note that Auslander algebras are useful when considering artin algebras of finite representation type, since there is a bijection between the Morita equivalence classes of such algebras  and the Morita equivalence classes of Auslander algebras (cf~\cite{ARS}).

The Hopf algebra structure on an algebra $\Lambda$ conveys an additional structure on the Grothendieck groups $K_0(\Lambda)$ and $\overline{K}_0(\Lambda) $ of isomorphism classes of projective (respectively all) indecomposable modules, since the tensor product over the base ring $k$ of two $\Lambda-$modules is again a $\Lambda-$module, via the comultiplication of $\Lambda.$ Furthermore, there is a one-to-one correspondance between the indecomposable modules over any algebra and the indecomposable projective modules over its Auslander algebra; in the case of a Hopf algebra, therefore, the Grothendieck group of projective modules of $\gl$ is endowed with a multiplicative structure. However, this correspondance does not preserve the underlying vector spaces, and this multiplicative structure doesn't seem natural.

In this paper, we study the example of the Taft algebras; they are Hopf algebras which are neither commutative, nor cocommutative. They are interesting for various reasons; for instance, $\Lambda _p$ is an example of a non-semisimple Hopf algebra whose dimension is the square of a prime (cf~\cite{M1}). They are of finite representation type; furthermore, when $n$ is odd, $\lan$ is isomorphic to the half-quantum group \uq  ($q$ primitive $n^{th}-$root of unity), and is the only half-quantum group $u_q^+(\frak{g})$ at a root of unity which is not of wild representation type (cf~\cite{C1}). Then, for each $n,$ $\lan$ is not braided, but its Grothendieck group is a commutative ring nonetheless (cf~\cite{C2}, \cite{G}). 

These examples show that the non-commutative, non-cocommutative Hopf algebra structure of $\lan$ does not yield a natural multiplicative structure on its cyclic homology. There is a product, however, obtained by transporting that of $K_0(\Lambda)$ via the Chern characters, which are onto.

The paper is organized as follows: first, we recall the Hochschild homology of truncated quiver algebras, among which are the Taft algebras, established by Sköldberg in~\cite{S}. Then, we compute the cyclic homology of these algebras, using the fact that they are graded. In the next section, we study the Auslander algebras of the Taft algebras: their quivers (which are the Auslander-Reiten quivers of the $\lan),$ their Hochschild homology, and their cyclic homology. Finally, we compute the Chern characters of the $\lan$ and of the $\gln.$

\section{The Taft algebras}
\subsection{Hochschild homology of truncated quiver algebras}

In this paragraph, $k$ is any commutative ring, except when a root of unity is needed.

Let  $\Delta_n$ be the following quiver (the $n-$crown): it has $n$ vertices $e_0, \ldots , e_{n-1},$ and  $n$ edges $a_0, \ldots , a_{n-1},$ each edge $a_i$ going from the vertex $e_i$ to the vertex $e_{i+1}$ for $0 \leq i \leq n-2,$ and the edge $e_{n-1}$ going from $e_{n-1}$ to $e_0,$ as follows:$$\epsfysize=4cm\epsfbox{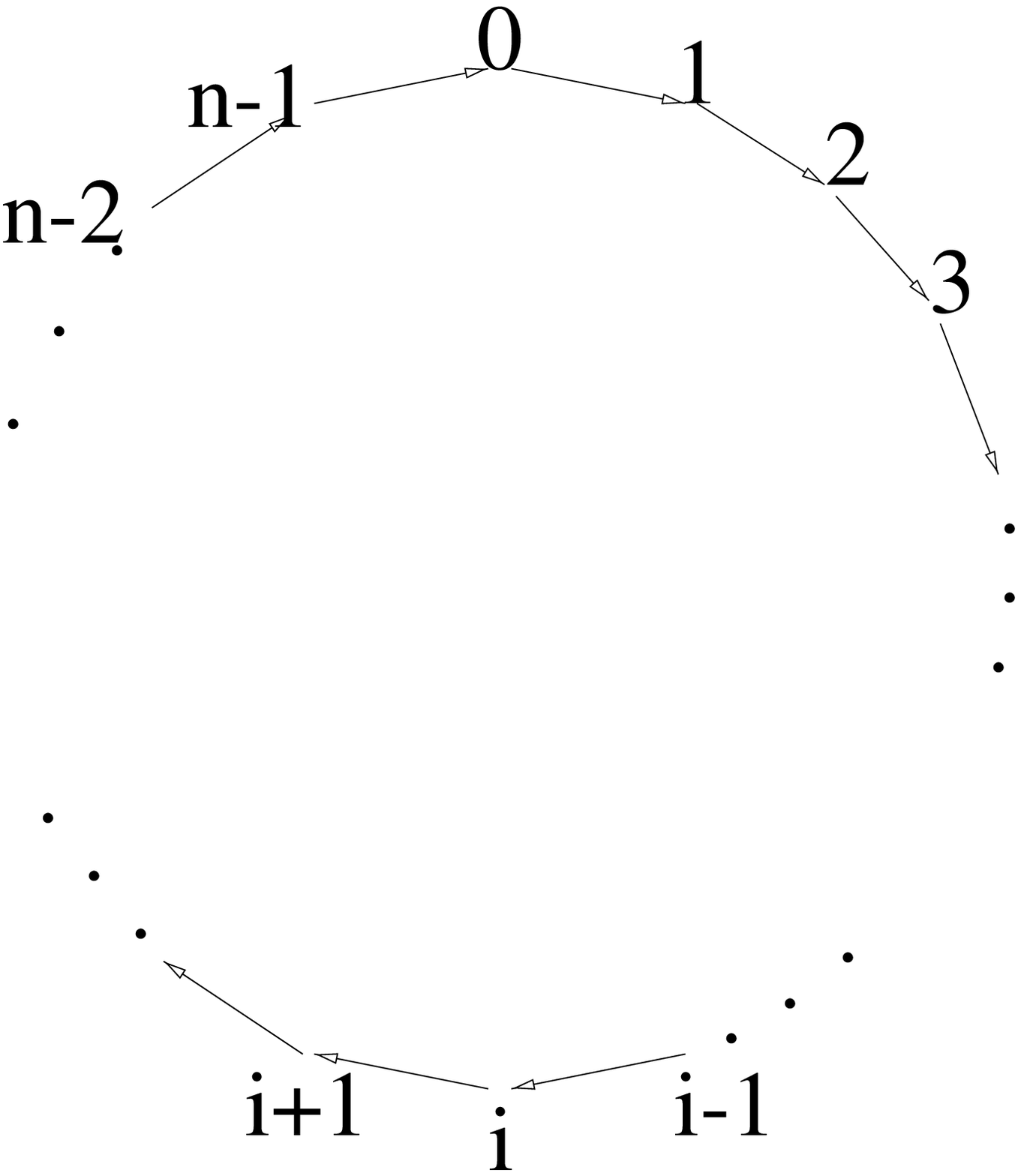}$$ Let $\m$ be the ideal in the path algebra of $\Delta_n$ generated by the paths of length 1.

\begin{dfn} The Taft algebra $\lan$ is the quotient of the path algebra $k\Delta_n$ by the ideal $\m^n.$ When $k$ contains a primitive $n^{th}$ root of unity $\zeta,$ it can also be described with generators and relations as follows: $\lan$ is the algebra generated by two elements $g$ and $x,$ subject to the relations $g^n=1, \  x^n=0,$ and $xg=\zeta gx.$ 
\end{dfn}

Still in the case where $k$ contains a  primitive $n^{th}$ root of unity $\zeta,$ the algebra $\Lambda _n =k \Delta / \m^n$ is a Hopf algebra (see~\cite{C1}), with the following structure maps:
\begin{eqnarray*} \varepsilon (e_i)=\delta _{i,0}, &&  \varepsilon (a_i)=0, \\
\Delta (e_i)= \sum_{j+k=i} e_j \ot e_k, && \Delta (a_i)= \sum_{j+k=i}( e_j \ot a_k+q^k a_j \ot e_k), \\
S(e_i)=e_{-i}, && S(a_i) = -q^{i+1} a_{-i-1}, 
\end{eqnarray*} where $\delta$ is the Kronecker symbol.

In \cite{S}, for any quiver $\Gamma,$ Sköldberg computes the Hochschild homology of the algebra $A:=k\Gamma/\m^n$ with coefficients in itself, when $n \geq 2;$ to state this result,  we shall need some notation: let $\mathcal{C}$ denote the set of cycles in the quiver $\Gamma,$ and for any cycle $\gamma$ in $\mathcal{C},$ let $L(\gamma)$ denote its length. There is a natural action of the cyclic group $\left< t_{\gamma} \right>$ of order $L(\gamma)$ on $\gamma;$ let $\overline{\gamma}$ denote the orbit of $\gamma$ under this action, and let  $\overline{\mathcal{C}}$ denote the set of orbits of cycles.

\begin{theo} [\cite{S}]\label{S} Set $q=cn+e,$ for $0 \leq e \leq n-1 \ (n \geq 2).$ Then:
$$HH_{p,q} (A) = \begin{cases} 
k^{a_q} & if \ 1 \leq e \leq n-1 \ and \ 2c \leq p \leq 2c+1, \\
\bigoplus _{r|q} (k^{(n\wedge r)-1 } \oplus Ker (. \frac{n}{n\wedge r} : k \rightarrow k))^{b_r} & if \ e=0, \ and \ 0<2c=p, \\
\bigoplus _{r|q} (k^{(n\wedge r)-1 } \oplus Coker (. \frac{n}{n\wedge r} : k \rightarrow k))^{b_r} & if \ e=0, \ and \ 0<2c-1=p, \\
k^{\# \Gamma_0} & if \ p=q=0, \ and \\
0 & otherwise,
\end{cases}  $$ where $a_q$ is the number of cycles of length $q$ in $\overline{\mathcal{C}},$   $b_r$ is the number of cycles of length $r$ in $\overline{\mathcal{C}}$ which are not powers of smaller cycles, and $n \wedge r$ is the greatest common divisor of $n$ and $r.$
\end{theo}

\begin{ex}\label{Hochtaft} The Hochschild homology of the Taft algebra $\Lambda _n$ is given by:
\begin{eqnarray*}
HH_{p,cn} (\Lambda _n )& =&k^{n-1} \mbox{ if $p=2c$ or $p=2c-1,$} \\
HH_{0,0} (\Lambda _n ) &=& k^n \\
HH_{p,q} (\Lambda _n ) &=& \ 0 \mbox{ in all other cases}.
\end{eqnarray*} \end{ex}


We shall now look at the cases $n=0$ and $n=1.$ 
The case $n=1$ is quite simple: $k\Gamma/\m$ is equal to $k \Gamma _0 \cong \times _{s \in \Gamma_0} ks,$ so that $$HH_p(k\Gamma /\m)=\bigoplus _{s \in \Gamma_0} HH_p(ks) = \begin{cases} \bigoplus _{s \in \Gamma_0} ks & \mbox{ if $ p=0,$} \\ 0 & \mbox{ if $p>0.$} \end{cases}$$

Finally, for the case $n=0,$ we can  state:

\begin{prop} The Hochschild homology of $k\Gamma$ is given by: $$\begin{cases} HH_0 (k\Gamma)=k\overline{\mathcal{C}}, \\ HH_1 (k\Gamma)= \{\sum_{i=0}^{L(\gamma) -1}t_{\gamma}^i (\gamma)/\ \gamma \in \overline{\mathcal{C}}, L(\gamma) \geq 1\} \\  HH_p (k\Gamma)=0 \mbox{ if $p \geq 2$}. \end{cases}$$
\end{prop}

\begin{proof} We shall use the following resolution (see for instance~\cite{C2}):
\begin{lemma}[\cite{C2} theorem 2.5] \label{resC} There is a $k\Gamma-$bimodule projective resolution of $k\Gamma$ given by $$\ldots  0\longrightarrow k\Gamma \ot _{k\Gamma_0} k\Gamma_1 \ot_{k\Gamma_0} k\Gamma \longrightarrow  k\Gamma \ot _{k\Gamma_0} k\Gamma \longrightarrow  k\Gamma  \longrightarrow 0. $$ \end{lemma}
Tensoring by $k\Gamma$ over $k\Gamma^e$ yields the following complex: $$\begin{array}{ccccc}
\ldots \longrightarrow 0 \longrightarrow &  k\Gamma \ot _{k\Gamma_0^e} k\Gamma_1& \stackrel{\delta}{\longrightarrow} &  k\mathcal{C} & \longrightarrow 0. \\
&\nu \ot a & \mapsto & \nu a - a \nu & \end{array}$$ 

The space $HH_0(k\Gamma)$ is generated by the cycles in $\mathcal{C},$ subjected to the relations given by the image of $\delta.$ Since $\delta (\nu \ot a)=\nu a - t_{\nu a }(\nu a ),$ the relations identify two cycles in the same orbit, and  $HH_0(k\Gamma)=k\overline{\mathcal{C}}.$

The complex is $\overline{\mathcal{C}}-$graded; therefore, to find $HH_1(k\Gamma)=\ker \delta,$ it is sufficient to consider elements of type $x=\sum_{i=0}^{L(\gamma)-1} \lambda_i t_{\gamma}^i (\gamma)$, in which $\nu \ot a$ is identified with $\nu a,$ and the $ \lambda_i$ belong to $k.$ We have $\delta (x) =0$ iff $ \lambda_0 = \lambda_1= \ldots = \lambda_{L(\gamma)-1},$ and the result follows. $\square$ \end{proof}

\subsection{Cyclic homology of graded algebras}

In this paragraph, $k$ is a commutative ring which contains $\Bbb{Q}.$ When $A$ is a graded $k-$algebra, Connes' SBI exact sequence splits in the following way:

\begin{theo} [\cite{L} Theorem 4.1.13]\label{SBIgr} Let $A$ be a unital graded algebra over $k$ containing $\Bbb{Q}.$ Define $\overline{HH}_p (A) = HH_p(A) / HH_p (A_0)$ and $\overline{HC}_p (A) = HC_p(A) / HC_p (A_0).$ Connes' exact sequence for $\overline{HC}$ reduces to the short exact sequences: $$0 \rightarrow \overline{HC}_{n-1} \rightarrow \overline{HH}_n \rightarrow \overline{HC}_n \rightarrow 0.$$  
\end{theo}

This will enable us to compute the  cyclic homology of truncated quiver algebras. Let us first consider the cases $n=0$ and $n=1.$ Combining the results for Hochschild homology and theorem~\ref{SBIgr} yields the following:

\begin{prop}\label{cycl01} The cyclic homologies of $k\Gamma$ and of $k\Gamma/\m$ are given by: \begin{eqnarray*} 
HC_{2c}(k\Gamma /\m)&=& \oplus_{s\in \Gamma_0} ks\\
HC_{2c+1}(k\Gamma /\m)&=& 0 \\
&\mathrm{and}\\
HC_{0}(k\Gamma)&=& k\overline{\mathcal{C}} \\
HC_{2c}(k\Gamma )&=& k^{\#\Gamma_0} \\
HC_{2c+1}(k\Gamma )&=& 0,
\end{eqnarray*} for all $c \in \Bbb{N}.$
\end{prop}

The case $n\geq 2$ involves the same methods:

\begin{prop} Suppose $n \geq 2.$ Then: 
\begin{eqnarray*}
\mathrm{dim}_k \,HC_{2c} (k\Gamma / \m ^n) &=& \#\Gamma _0 + \sum_{e=1}^{n-1} a_{cn+e} - \sum_{\tiny \begin{array}{c} r|(c+1)n \\ 
n \notin r \Bbb{N} \end{array}} (r \wedge n -1)b_r \\ 
\mathrm{dim}_k \,HC_{2c+1} (k\Gamma / \m ^n) &=&  \sum_{r|n} ( r -1) b_r.
\end{eqnarray*}
\end{prop}

\begin{proof} In the first place, $A_0$ is equal to $k\Gamma_0,$ so that we know the homologies of  $A_0$ (see proposition \ref{cycl01}). Next, we have $HC_0(A) = HH_0(A) = k^{\#\Gamma_0 + \sum_{e=1}^{n-1} a_e}.$ Then, using theorem \ref{SBIgr}, we get the following formula: $$dim_k HC_{2c}(A) +dim_k HC_{2c+1}(A)=\#\Gamma_0 + \sum_{r|(c+1)n}(r \wedge n -1)b_r + \sum_{e=1}^{n-1} a_{cn+e}.$$ In particular, $dim_k HC_{1}(A) = \sum_{r|n}(r-1)b_r.$

An induction on $c$ yields the result. $\square$ \end{proof}

\begin{cor}\label{cycltaft} When $\Gamma =\Delta_n$ is the $n-$crown, then the results are:
$$\begin{cases} HC_{2c} (\Lambda _n ) = k^n, \\
HC_{2c+1} (\Lambda _n ) = k^{n-1} \mbox{ for $c \in \Bbb{N}.$}
\end{cases}$$
\end{cor}

\begin{rmk} The cyclic cohomology of these truncated quiver algebras has been computed in \cite{BLM} and \cite{Li}.
\end{rmk}

\section{The Auslander algebra of $\Lambda _n$}

In this section, $k$ is an algebraically closed field.

\subsection{The quiver of the Auslander algebra}

\begin{dfn} Let  $\Lambda $ be a finite-dimensional basic algebra over an algebraically closed field $k,$ with only a finite number of isomorphism classes of indecomposable modules. The Auslander algebra of  $\Lambda$ is: $$\gl := End _{ \Lambda} (\oplus_{M \in \mathrm{ind}M} M)^{op},$$ where $\mathrm{ind}M$ is the set of isomorphism classes of indecomposable $\Lambda-$modules.
 \end{dfn}

The algebra $\gl$ has a quiver, which is the opposite of the Auslander-Reiten quiver of $\Lambda.$ The relations on this quiver are given by the mesh relations (see~\cite{ARS} p232).

When the algebra $\Lambda$ is $\Lambda _n,$ with $n \geq 1,$ the quiver is the following:

$$\epsfysize=9cm\epsfbox{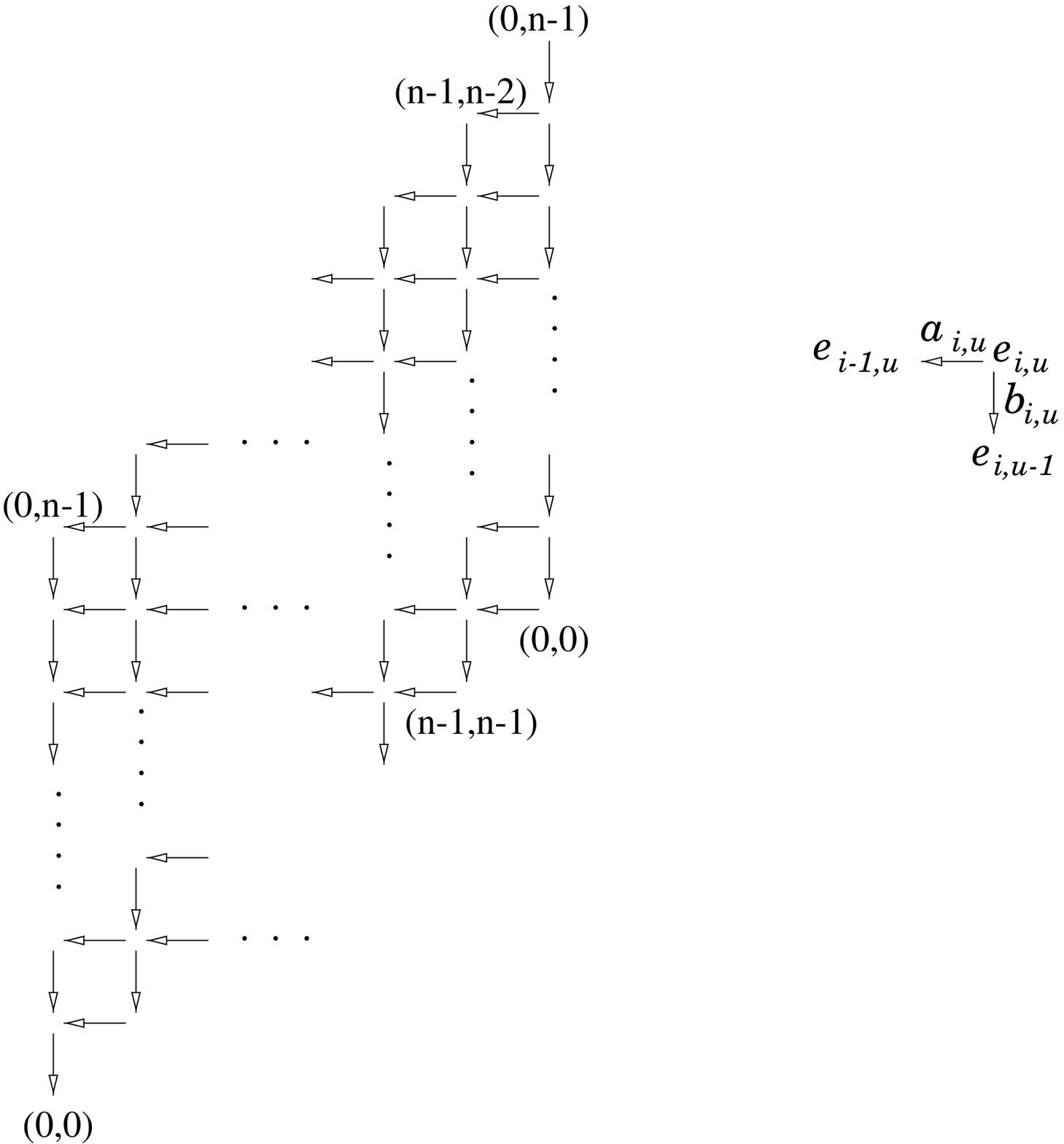}$$
where both vertical outer edges are identified (the quiver is on a cylinder: see~\cite{GR}). Let $Q$ denote this quiver,  let $\{e_{i,u}/(i,u) \in \Bbb{Z}/n\Bbb{Z} \times \Bbb{Z}/n\Bbb{Z}\}$ be the set of vertices of $Q,$ and let $\{a_{i,u}\, ; \,b_{i,u}/ \,(i,u) \in \Bbb{Z}/n\Bbb{Z} \times \Bbb{Z}/n\Bbb{Z}\}$ be the set of edges of $Q,$ as in the figure above. 

The mesh relations on this quiver are: $a_{i,i-2}b_{i,i-1}=0$ for all $i \in \Bbb{Z}/n\Bbb{Z}$ (the composition of two edges of any `triangle' under the top diagonal is zero), and $a_{i,iu-1}b_{i,u}+b_{i-1,u}a_{i,u}=0$ for all $i$ and $u$ in $\Bbb{Z}/n\Bbb{Z}$ (the squares are anticommutative). 

The algebra $\gln$ is the quotient of the path algebra $kQ$ by the ideal generated by these relations.

\begin{rmk} The algebra $\gln$ is not a Hopf algebra, since its quiver is not a Cayley graph (see~\cite{GS} theorem 2.3; in relation to this, see also  \cite{CR}, in which the authors study the case without relations). Another argument can be given: if $\gln$ were a Hopf algebra, it would be selfinjective as an algebra (cf~\cite{M2}), therefore of homological dimension 0 or $\infty;$  however, the homological dimension of an Auslander algebra is at most 2 in general (see~\cite{ARS}), and in this case it is not zero, since $\lan$ itself is not semisimple (\cite{ARS} proposition 5.2 p211).
\end{rmk}


\subsection{Projective resolutions of simple $\gln-$modules}

Let $P_{i,u}$ denote the indecomposable projective $\gln-$module at the vertex $e_{i,u},$ and let $S_{i,u} = top (P_{i,u})$ be the corresponding simple module. These modules are described on the quiver by:

$$\epsfysize=5cm\epsfbox{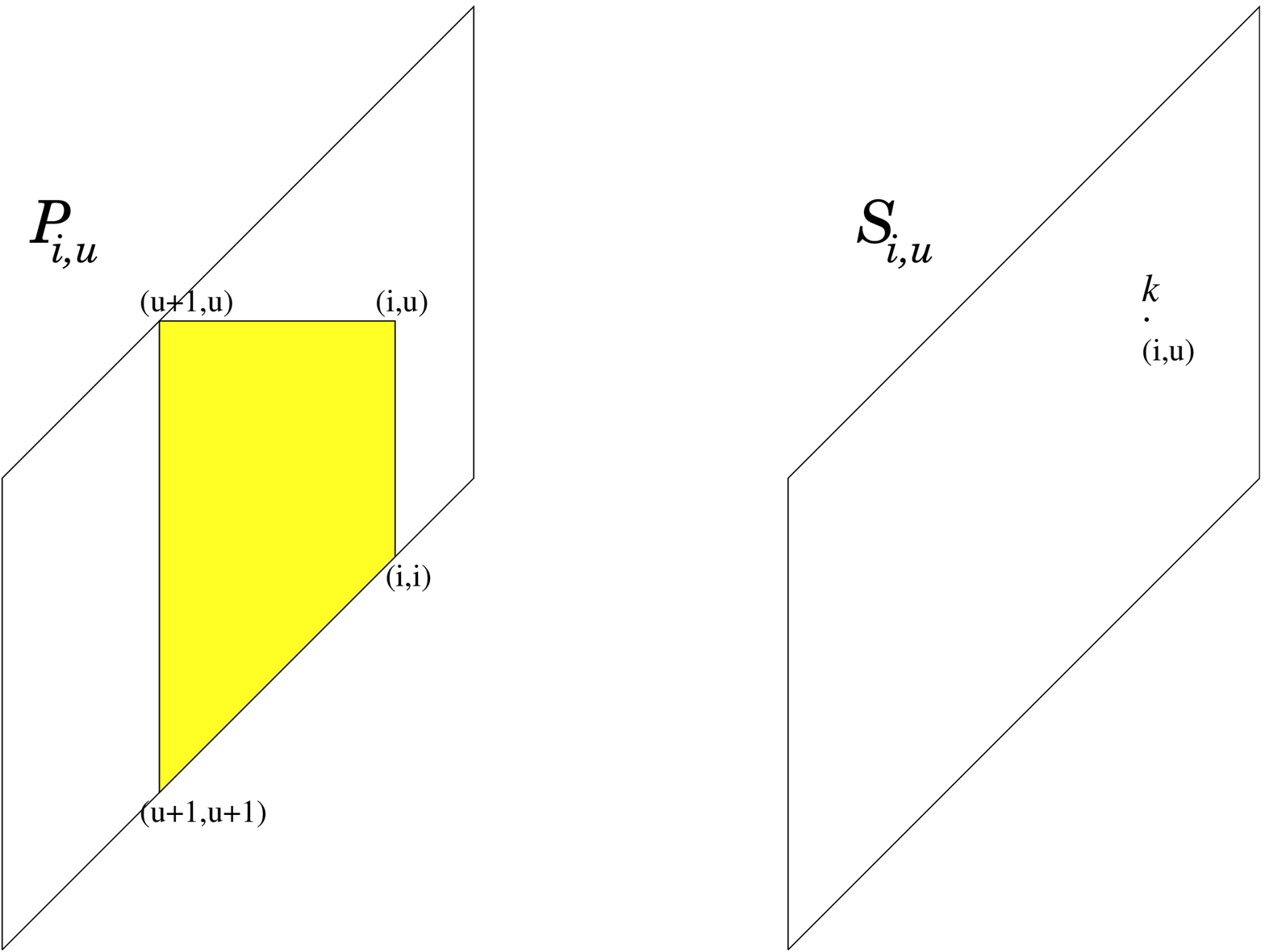}$$ where each vertex in the shaded part is represented by $k,$ and each edge in the shaded part is represented by $id.$ Everywhere else, the vertices and edges are represented by 0 (see also~\cite{GR} section 2).

We can  compute the minimal projective resolutions of the simple $\gln-$modules; the results are as follows:

\begin{prop}\label{sple-res} The projective resolutions of the simple modules which are obtained by successive projective covers are: $$\begin{array}{ccccccccr} 
&&0 \ \longrightarrow & P_{i-1,i} &  \longrightarrow & P_{i,i} & \longrightarrow & S_{i,i} &  \longrightarrow  0 \\
0  \longrightarrow & P_{i-1,i-2} & \longrightarrow & P_{i,i-2} & \longrightarrow & P_{i,i-1} & \longrightarrow & S_{i,i-1} & \longrightarrow 0 \\
0  \longrightarrow & P_{i-1,i-j-1} & \longrightarrow & P_{i-1,i-j} \oplus P_{i,i-j-1} & \longrightarrow & P_{i,i-j} & \longrightarrow & S_{i,i-j} & \longrightarrow 0
\end{array}$$ for $2 \leq j \leq n-1.$ 
\end{prop}

\begin{proof} We consider only the $S_{n-1,u},$ because the other cases may be obtained by translating the quiver along the cylinder on which it lies. We then determine the radical of $P_{n-1,u},$ and the latter's projective cover, through their representations on the quiver. Iterating this process until the radical obtained is projective yields the result.
 \end{proof}

In particular, this enables us to compute the $Ext^p_{\gln} (S\, ;T),$ where $S$ and $T$ are two simple modules:

\begin{prop} Let $S$ be a simple $\gln-$module. Then:
\begin{eqnarray*}
Ext^1_{\gln} (S_{i,u}\, ; S)& = &\begin{cases} k & \mbox{ if $S=S_{i,u}$,} \\ 0 & \mbox{ if $S \neq S_{i,u}$,}\end{cases} \\
Ext^2_{\gln} (S_{i,i}\, ; S)& =& 0\\
Ext^2_{\gln} (S_{i,i-j}\, ; S)& =&\begin{cases} k & \mbox{ if $S=S_{i-1,i-j-1}$,} \\ 0 & \mbox{ if $S \neq S_{i-1,i-j-1}$,}\end{cases} \mbox{ for 1 $\leq j \leq n-1$} \\
Ext^p_{\gln} (S_{i,u}\, ; S)& = & 0 \mbox{ for $ p \geq 3.$}
\end{eqnarray*}
\end{prop}

\begin{proof} The $Ext^1$ are known, owing to Schur's lemma. The other $Ext^i$ are obtained from the resolutions of proposition~\ref{sple-res}, applying the functor $Hom_{\gln}(-; S).$ 
 \end{proof}


\subsection{Hochschild and cyclic homologies of $\gln$}

We are going to use the previous results to compute a minimal projective resolution of $\gln$ as a $\gln-$bimodule, due to Happel (he does it in the general situation in \cite{H} 1.5.):

\begin{theo}[\cite{H}]\label{happel} If $$\ldots \rightarrow R_p \rightarrow R_{p-1} \rightarrow \ldots \rightarrow R_1 \rightarrow R_0 \rightarrow \gln \rightarrow 0$$ is a minimal projective resolution of $\gln$ as a $\gln-$bimodule, then $$R_p= \bigoplus_{\tiny \begin{array}{c} (i,u) \\(j,v) \end{array}} (\gln e_{j,v} \ot e_{i,u} \gln )^{\mathrm{dim} _k Ext_{\gln}^p (S_{i,u} \, ; S_{j,v})}.$$

Specifically: \begin{eqnarray*} 
R_0 &=& \bigoplus_{\tiny (i,u)} \gln e_{i,u} \ot e_{i,u} \gln \\
R_1 &=& \bigoplus_{\tiny (i,u)} [ (\gln e_{i-1,u} \ot e_{i,u} \gln) \oplus  (\gln e_{i,u-1} \ot e_{i,u} \gln) ] \\
R_2 &=& \bigoplus_{\tiny \begin{cases}(i,u) \\ i \neq u \end{cases}} \gln e_{i-1,u-1} \ot e_{i,u} \gln \\
R_p & =& 0 \mbox{ if $p \geq 3.$}
\end{eqnarray*}
\end{theo}

Applying the functor $\gln \ot _{\tiny \gln - \gln} -,$ we obtain a complex: $$\ldots 0 \longrightarrow \ldots \longrightarrow 0 \longrightarrow k Q_0 \longrightarrow 0.$$ This yields:

\begin{prop} The Hochschild homology of $\gln$ is: $$\begin{cases} HH_0 (\gln) = kQ _0 \cong k^{n^2} \\
HH_p (\gln) =0 & \forall p \in \Bbb{N}^*, \end{cases}$$ and hence the cyclic homology of $\gln$ is: $$\begin{cases} HC_{2p} (\gln) = kQ _0 \cong k^{n^2} \\
HH_{2p+1} (\gln) =0  \end{cases}  \forall p \in \Bbb{N}.$$
\end{prop}

\begin{rmk} There doesn't seem to be any connection between these results and those for $\lan$ (see the example on page \pageref{Hochtaft} and corollary \ref{cycltaft}).
\end{rmk}


\section{Chern characters of $\lan$ and $\gln$}

Let $K_0(\lan)$ (resp. $ K_0(\gln)$) be the Grothendieck group of projective $\lan -$modules (resp. $\gln -$modules). We are interested in the Chern characters $ch_{0,p} : K_0(\lan) \rightarrow  HC_{2p} (\lan )$ (resp. $ K_0(\gln) \rightarrow HC_{2p} (\gln)$). We shall write $[P_j]$ (resp. $[P_{i,u}]$) for the isomorphism class of the projective module at the vertex $e_j$ (resp. $e_{i,u}$).

Set $\sigma ^p = (y_p,z_p,\ldots, y_1,z_1,y_0) \in \Bbb{N}^{2p+1}$ with $y_p=(-1)^p (2p)!/p!$ and $z_p=(-1)^{p-1} (2p)!/2(p!).$ There is a system of generators of $HC_{2p} (\lan)$ (resp. $HC_{2p} (\gln)$) given by the following set: $$\{\sigma _i^p := \sigma ^p (e_i, \ldots , e_i)\in (Tot CC(\lan))_{2p}\; /\;i=0, \ldots, n-1 \}$$ (resp. by $\{\sigma _{i,u}^p := \sigma ^p (e_{i,u} , \ldots , e_{i,u})\in (Tot CC(\gln))_{2p}\; /\;i, u \in \{0,1, \ldots, n-1\}  \}).$

Consider the elements \begin{eqnarray*}\epsilon_j : \lan & \longrightarrow & \lan \\ \lambda & \mapsto & \lambda e_j \end{eqnarray*} in $\mathcal{M}_1(\lan)$ and \begin{eqnarray*}\epsilon_{i,u} : \gln & \longrightarrow & \gln \\ \lambda & \mapsto & \lambda e_{i,u} \end{eqnarray*} in $\mathcal{M}_1(\gln);$ their ranges are the corresponding projective modules. Then by definition of the Chern characters (see \cite{L} 8.3.4), we have: $$\begin{array}{l} 
ch_{0,p} ([P_j]) = ch_{0,p} ([\epsilon_j]) := \mathrm{tr} (c(\epsilon_j))= \sigma _j^p  \ \mbox{in $HC_{2p}(\lan)$}\\
ch_{0,p} ([P_{i,u}])= ch_{0,p} ([\epsilon_{i,u}]) = \sigma _{i,u}^p.
\end{array}$$ using the isomorphisms $\mathcal{M}_m(\Lambda) \cong \mathcal{M}_m(k) \ot \Lambda.$ Here, $$c(\epsilon_j)=(y_p \epsilon_j ^{\ot 2p+1} , z_p \epsilon_j ^{\ot 2p}, \ldots , z_1 \epsilon_j ^{\ot 2}, y_0 \epsilon_j) \in \mathcal{M}(\gln)^{\ot 2p+1} \oplus \ldots \oplus \mathcal{M}(\gln).$$
 
\begin{rmk} There is a decomposition formula for the tensor product of indecomposable modules on $\lan$ (see~\cite{C2}, \cite{G}). From this formula, we get inductively: $$ch_{0,p} ([L_1] \ot \ldots \ot [L_r]) = \left( 1/n^2\right) \prod_{i=1}^r \left( \mathrm{dim} L_i \right) \; (\sigma _0^p ,\ldots , \sigma _{n-1} ^p), \ \mbox{for $r\geq 2,$}$$ where the $L_i$ are arbitrary projective $\lan-$modules. Unfortunately, this product in the cyclic homology doesn't seem natural.
\end{rmk}

\begin{rmk} Let $\overline{K}_0 (\lan)$ be the Grothendieck group of all $\lan-$modules (not just the projective ones). Then $\overline{K}_0 (\lan) \cong K_0(\gln).$ Hence, if $N_{i,u} $ is the indecomposable $\lan-$module which starts at the vertex $i$ and ends at the vertex $u,$ it corresponds to the projective $\gln-$module $P_{i,u},$ and we get a map: 
\begin{eqnarray*} 
\overline{K}_0 (\lan)& \longrightarrow & HC_{2p} ( \gln) \\ 
N_{i,u} & \mapsto & \sigma_{i,u}^p. 
\end{eqnarray*}
\end{rmk}

\begin{rmk} Although $\gln$ is not a Hopf algebra, its Grothendieck group $K_0(\gln)$ does have a ring structure, albeit not natural: for every $[P]$ in $K_0(\gln),$ there exists a $[B]$ in $K_0(\lan)$ such that $[P]=[Hom_{\lan}(M,B)],$ where $M$ is the sum of all isomorphism classes of indecomposable $\lan-$modules. If $[Q]=[Hom_{\lan}(M,C)]$ is another element in  $K_0(\gln),$ we can set $$[P].[Q]= [Hom_{\lan}(M,B\ot _k C)]$$ (the vector space $B\ot _k C$ is a $\lan-$module since $\lan$ is a Hopf algebra).  In fact, using the decomposition in \cite{C2} and \cite{G}, the product can be written: $$[P_{i,u}][P_{j,v}]=  \begin{cases} \sum_{l=0}^{v-j} [P_{i+j+l, u+v -l}] & \mbox{ if $u+v-(i+j) \leq n-1$}\\ \sum_{l=0}^{e} [P_{i+j+l, u+v +l-1}] + \sum_{m=e+1}^{v-j} [P_{i+j+m, u+v -m}] & \mbox{ if $e:=u+v-(i+j) -( n-1) \geq 0.$}
\end{cases}$$
\end{rmk}



\begin{thebibliography}{99}

\bibitem[ARS]{ARS} \textsc{Auslander, M.,  Reiten, I. and Smalø, S.O.,} Representation theory of Artin algebras, \textit{Cambridge Studies in Advanced Mathematics, 36. Cambridge University Press, Cambridge,} (1997). \\

\vspace*{-15pt}

\bibitem[BLM]{BLM} \textsc{Bardzell, M.J., Locateli, A.C., Marcos, E.N.,} On the Hochschild cohomology of truncated cycle algebras, \textit{Comm. Algebra} \textbf{28} (2000), no. 3, pp 1615-1639.\\

\vspace*{-15pt}

\bibitem[C1]{C1}  \textsc{Cibils, C.,} Half-quantum groups at roots of unity, path algebras, and representation type. Internat. Math. Res. Notices
1997, no. 12, 541--553. \\

\vspace*{-15pt}

\bibitem[C2]{C2} \textsc{Cibils, C.,} A quiver quantum group, \textit{Comm. Math. Phys.} \textbf{157} (1993), no. 3, pp 459-477. \\

\vspace*{-15pt}

\bibitem[CR]{CR} \textsc{Cibils, C. and Rosso, M.,} Alg\`ebre des chemins quantiques,  \textit{Adv. Math.} \textbf{125} (1997), pp 171-199. \\

\vspace*{-15pt}

\bibitem[GR]{GR} \textsc{Gabriel, P. and  Riedtmann, Ch.,} Group representations without groups,  \textit{Comment. Math. Helv.}  \textbf{54} (1979), no. 2, pp 240-287. \\

\vspace*{-15pt}

\bibitem[GS]{GS} \textsc{Green, E.L. and Solberg,  Ø.,} Basic Hopf algebras and quantum groups, \textit{Math. Zeit.} \textbf{299} (1998), no. 1, pp 45-76. \\ 

\vspace*{-15pt}
 
\bibitem[G]{G} \textsc{Gunnlaugsd\'ottir, E.,} Monoidal structure of the category of $u_q^+-$modules,   \textit{preprint no. 00-14} (2000) \textit{of the Math. Department in Montpellier}.\\

\vspace*{-15pt}

\bibitem[H]{H} \textsc{Happel, D.,} Hochschild cohomology of finite-dimensional algebras, \textit{Séminaire d'Algèbre Paul Dubreil et Marie-Paul Malliavin, 39ème Année (Paris, 1987/1988),} pp 108-126, \textit{Lecture Notes in Math.} \textbf{1404}, \textit{Springer, Berlin-New York} (1989).\\

\vspace*{-15pt}

\bibitem[Li]{Li} \textsc{Locateli, A. C.,} Hochschild cohomology of truncated quiver algebras, Comm. Algebra \textbf{27} (1999), no. 2, pp 645-664. \\

\vspace*{-15pt}

\bibitem[L]{L} \textsc{Loday, J-L.,} Cyclic homology, Appendix E by Mar\'\i a O. Ronco, \textit{Springer-Verlag, Berlin,} (1992). \\

\vspace*{-15pt}

\bibitem[M1]{M1} \textsc{Montgomery, S.,} Classifying finite-dimensional semisimple Hopf algebras,  \textit{Contemp. Math.} \textbf{229} (1998), pp 265-279, \textit{Amer. Math. Soc., Providence, RI.} \\

\vspace*{-15pt}

\bibitem[M2]{M2} \textsc{Montgomery, S.,} Hopf algebras and their actions on rings,  \textit{CBMS Regional Conference Series in Mathematics} \textbf{82}, \textit{Published for the Conference Board of the Mathematical Sciences, Washington, DC; by the American Mathematical Society, Providence, RI} (1993). \\

\vspace*{-15pt}

\bibitem[S]{S}  \textsc{Sköldberg, E.,} The Hochschild homology of truncated and quadratic monomial algebras, \textit{J. London Math. Soc. (2)} \textbf{59} (1999), no. 1, pp 76-86 (1999).

\end{thebibliography}
\end{document}